\documentclass[12pt]{amsart}
\usepackage{epsfig}
\numberwithin{equation}{section}

\newtheorem{theorem}[equation]{Theorem}
\newtheorem{lemma}[equation]{Lemma}
\newtheorem{corollary}[equation]{Corollary}
\newtheorem{proposition}[equation]{Proposition}

\theoremstyle{definition}

\theoremstyle{definition}           

\theoremstyle{remark}
\newtheorem{remark}{Remark}

\newtheorem{claim}{Claim}

\newtheorem{assertion}{Assertion}

\newcommand{\eps}{\varepsilon}

\newcommand\vphi{\varphi}

\newcommand\al{\alpha}

\newcommand\be{\beta}
\newcommand\Si{\Sigma}
\newcommand\ga{\gamma}
\newcommand\Ga{\Gamma}

\newcommand\De{\Delta}

\newcommand\cN{\mathcal  N}

\newcommand\cS{\mathcal  S}

\newcommand\cK{\mathcal  K}

\newcommand\cF{\mathcal F}

\newcommand\bfC{\mbox {\bf  C}}

\newcommand\bfP{\mbox {\bf  P}}

\newcommand\bfQ{\mbox {\bf  Q}}

\newcommand\bfZ{\mbox {\bf  Z}}

\newcommand\nl{\newline}

\newcommand\wtl{\widetilde}

\newcommand\order{\mbox{\text{order}\/}}

\newcommand\GL{\text{GL}\/}
\newcommand\PGL{\text{PGL}\/}

\newcommand\id{\text{id}}

\newcommand\coeff{\text{coeff}\/}
\newcommand\dual{\text{dual}\/}

\newcommand\inv{^{-1}}
\def\mapright#1{\smash{\mathop{\longrightarrow}\limits^{{#1}}}}

\def\mapdown#1{\Big\downarrow\rlap{$\vcenter{\hbox{$#1$}}$}}

\def\inv{^{-1}}

\begin{document}
\title
{ Elliptic curves from sextics}
\author{Mutsuo Oka }
\address{
Department of Mathematics,
Tokyo Metropolitan University\newline\indent
Minami-Ohsawa, Hachioji-shi
 Tokyo 192-03, Japan}
\email{oka@comp.metro-u.ac.jp}
\date{November, 1999, first version}
\begin{abstract}
 Let $\cN$ be the moduli space of sextics
 with 3 (3,4)-cusps. The quotient moduli space $\cN/G$ is
one-dimensional and consists of  two components, 
$\cN_{torus}/G$
and 
$\cN_{gen}/G$. 
By quadratic transformations,
they  are  transformed
into 
 one-parameter families  $C_s$
and 
$D_s$ of cubic curves respectively.
First we study the geometry of $\cN_\eps/G,\eps=torus,gen$ and 
 their structure of elliptic fibration. Then  we
study   the Mordell-Weil torsion
groups   of  cubic curves $C_s$ over $\bfQ$ and $D_s$ over $\bfQ(\sqrt{-3})$
respectively.
We show that $C_{s}$  has  the torsion group 
$\bfZ/3\bfZ$ for a generic $s\in \bfQ$  and it also 
contains
subfamilies  which coincide  with the universal families given by 
Kubert \cite{Kubert} with the torsion groups
 $\bfZ/6\bfZ$, $\bfZ/6\bfZ+\bfZ/2\bfZ$,
$\bfZ/9\bfZ$ or
$\bfZ/12\bfZ$. The cubic
curves $D_s$ has torsion 
$\bfZ/3\bfZ+\bfZ/3\bfZ$  generically but also 
$\bfZ/3\bfZ+\bfZ/6\bfZ$ for a subfamily  which is  parametrized by
$ \bfQ(\sqrt{-3}) $.
\end{abstract}
\maketitle 
\vspace{1cm}
\section{Introduction}
Let $\cN_3$ be the moduli space of 
sextics with 3 (3,4)-cusps as in \cite{dual}. For brevity, we
denote $\cN_3$ by $\cN$.  A sextic $C$ is called {\em of a
torus type} if its defining polynomial $f$ has the
expression
$f(x,y)=f_2(x,y)^3+f_3(x,y)^2$ for some polynomials 
$f_2,f_3$ of degree 2, 3 respectively. 
We denote by  $\cN_{torus}$   the component
of $\cN$ which consists of curves of a torus type and by
$\cN_{gen}$ the curves of a general type (=not  of a torus type).
We denote  the dual curve of $C$ by $C^*$.
 Let $G=\PGL(3,\bfC)$. The quotient moduli space is
by definition the
quotient space of the moduli space  by the action of $G$.

In \S 2, we study the quotient moduli space $\cN/G$. We will show that
 $\cN/G$ is one dimensional
and it has  two components
$\cN_{torus}/G$ and  $\cN_{gen}/G$ which
consist  of sextics of  a torus type and sextics of a general
type respectively. After giving normal forms of these components
 $C_s,s\in \bfP^1(\bfC)$
and 
$D_s, s\in \bfP^1(\bfC)$, 
we show that the family {\em
 $ C_s$
 contains a  unique sextic  $C_{54}$ which is self dual}
(Theorem \ref{uniqueness})  and  
{\em $C_{54}$ has 
an involution which is associated with the Gauss map
 }(Proposition \ref{involution6}).

 In section 3, we study the structure of the elliptic
fibrations on the components $\cN_{\eps}/G,~\eps=torus, gen$
which are represented by the normal families
$C_s,s\in \bfP^1(\bfC)$
and 
$D_s, s\in \bfP^1(\bfC)$.
 Using 
 a quadratic  transformation
we  write these   families by  smooth cubic
curves
$C_s$ and $D_s$ which are defined by the following cubic polynomials.
\begin{eqnarray*} &C_s:~x^3-\frac1 4 s (x-1)^2+s y^2=0\qquad\qquad\\
& D_s:~ -8 x^3+1+(s+35)y^2-6 x^2+3x-6\sqrt{-3}y-3  
\sqrt{-3}x  \\&
\qquad\qquad\qquad -6\sqrt{-3}x^2-12 \sqrt{-3} xy+(s-35)xy=0
\end{eqnarray*}
 We show that $C_s,~s\in \bf\bfP^1(\bfC)$ (respectively 
 $D_s,~s\in \bf\bfP^1(\bfC)$) {\em has the structure
of rational elliptic surfaces}
$X_{431}$ (resp. $X_{3333}$) in the notation of \cite{Miranda}.

In section 4, we study their torsion subgroups of the Mordell-Weil group
of the
cubic families $C_s$ and $D_s$. The family $C_s$ is defined over $\bfQ$ and
$D_s$ is defined over
quadratic number field $\bfQ(\sqrt{-3})$. 
Both families enjoy   beautiful arithmetic properties. 
We will show that {\em the
torsion group $(C_s)_{tor}(\bfQ)$  is isomorphic to 
$ \bfZ/3\bfZ$ for a generic $s\in \bfQ$ but it has
subfamilies
$C_{\vphi_6(u)}$,  $C_{\vphi_{6,2}(r)}$, 
$C_{\vphi_9(t)}$ and $C_{\vphi_{12}(\nu)}$, $u,r,t,\nu\in \bfQ$
for which the Mordell-Weil torsion group are 
$\bfZ/6\bfZ,~\bfZ/6\bfZ+\bfZ/2\bfZ, ~\bfZ/9\bfZ$  and 
$\bfZ/12\bfZ$ respectively}.  Each of these  groups is parametrized 
by a rational function with $\bfQ$ coefficients which is defined over
$\bfQ$  and this parametrization  coincides, up to a linear fractional
change of
parameter, to the universal family given by Kubert in \cite{Kubert}.

As for $(D_s)_{tor}(\bfQ(\sqrt{-3}))$, we show that 
{\em $(D_s)_{tor}(\bfQ(\sqrt{-3})$ is generically isomorphic to 
$\bfZ/3\bfZ+\bfZ/3\bfZ$
but it also takes  $\bfZ/3\bfZ+\bfZ/6\bfZ$ 
for a subfamily $D_{\xi_6(t)}$ parametrized by
a rational function with coefficients in $\bfQ$ and 
defined on $\bfQ(\sqrt{-3})$}.

\vspace{.3cm}
\section{Normal forms of the moduli $\cN$}
We
consider the submoduli $\cN^{(1)}$ of the sextics whose
  cusps are at 
$O:=(0,0),A:=(1,1)$ and $B:=(1,-1)$.
As every sextic in $\cN$ can be represented by a curve
in $\cN^{(1)}$ by the action of $G$, we have $\cN/G\cong
\cN^{(1)}/G^{(1)}$ where  $G^{(1)}$ is the stabilizer
of $\cN^{(1)}$: $G^{(1)}:=\{g\in 
G;g(\cN^{(1)})=\cN^{(1)}\}$. By an easy computation, 
we   see that
$G^{(1)}$ is the semi-direct product of the group
$G_0^{(1)}$ and a finite group $\cK$,
isomorphic to the permutation group
$\cS_3$ where   $G_0^{(1)}$ is 
defined by 
\[
G_0^{(1)}:=\big\{ M=\left (
\begin{matrix}
a_1&a_2&0\\a_2&a_1&0\\ a_1-a_3&a_2&a_3\end{matrix}
\right )\in G; a_3(a_1^2-a_2^2)\ne 0\big\}
\] 

Note that $G_0^{(1)}$ is normal in $G^{(1)}$ and 
$g\in G_0^{(1)}$  fixes singular points
pointwise. The isomorphism
$\cK\cong \cS_3$ is given by identifying $g\in \cK$
as the  permutation of three singular locus $O,A,B$.
   We will study the normal forms of the quotient
moduli
$\cN/G\cong
\cN^{(1)}/G^{(1)}$.
\begin{lemma}\label{iso1}
For a given line $L:=\{y=bx\}$ with $b^2-1\ne 0$, 
there exists $M\in G_0^{(1)}$
such that $L^M$ is given by $x=0$.
\end{lemma}
{\em Proof.} By an easy computation, the image  of $L$
by the action of $M\inv$, where $M$ is as above,  is defined by
$(a_1-ba_2)y+(a_2-ba_1)x=0$. Thus we take $a_1=ba_2$. Then
$a_1^2-a_2^2=a_2^2 (b^2-1)\ne 0$ by the assumption. \qed

\begin{lemma} \label{iso2}The tangent
cone at $O$
 is not $y \pm x=0$ for  $C\in \cN^{(1)}$.
\end{lemma} 
{\em Proof.} Assume for
example that  $y-x=0$ is the
tangent cone of $C$ at $O$.
The  intersection multiplicity of
the line
$L_1:=\{y-x=0\}$ 
 and $C$  at $O$ is 4 and thus $L_1\cdot C\ge 7$, an obvious
contradiction to Bezout theorem. \qed

Let
$\cN^{(2)}$ be the subspace of $\cN^{(1)}$ consisting
of curves  $C\in \cN^{(1)}$ whose
 tangent cone at $O$ is given by
$x=0$.
Let $G^{(2)}$ be the stabilizer of $\cN^{(2)}$.
By Lemma
\ref{iso1} and Lemma \ref{iso2}, we have the isomorphism :
\begin{corollary}  $\cN^{(1)}/G^{(1)}\cong
\cN^{(2)}/G^{(2)}$.\end{corollary}

It is easy to see that $G^{(2)}$ is generated by the
group $G_0^{(2)}:=G^{(2)}\cap G_0^{(1)}$
 and an element 
$\tau$ of order two which is  defined by  $\tau(x,y)= (x,-y)$.
Note that 
\[G_0^{(2)}=\{
M=\left (
\begin{matrix}
a_1&0&0\\0&a_1&0\\ a_1-a_3&0&a_3\end{matrix}
\right)\in G_0^{(1)};\quad a_1 a_3\ne 0\}
\]

For $C\in \cN^{(2)}$, we associate complex numbers
$b(C),c(C)\in \bfC$ which are  the directions of the tangent cones of
$C$ at
$A,B$ respectively. This implies that  the lines
$y-1=b(C)(x-1)$ and 
$y+1 =c(C)(x-1)$ are  the tangent cones of 
$C$ at $A$ and $B$ respectively. We have shown that
$C\in \cN_{torus}^{(2)} $ if and only if $b(C)+c(C)=0$
and otherwise
   $C$ is  of a general type and  they satisfy
$c(C)^2+3c(C)-b(C)c(C)+3-3b(C)+b(C)^2=0$
(\S 4, \cite{dual}).

We consider the subspaces:
\[\cN_{torus}^{(3)}:=
\{C\in\cN_{torus}^{(2)};
b(C)=0\},\quad
 \cN_{gen}^{(3)}:=
\{C\in\cN_{gen}^{(2)};
b(C)=c(C)=\sqrt{-3}\}\]
and we put  $\cN^{(3)}:=\cN_{torus}^{(3)}\cup
\cN_{gen}^{(3)}$.
\begin{remark}  The common solution of  both
equations: $b+c=c^2+3c-bc+3-3b+b^2=0$ is $(b,c)=(1,-1)$ and
 in this
case,
$C$ degenerates into two non-reduced lines $(y^2-x^2)^2=0$ and a
conic. \end{remark}
\begin{lemma}
Assume that $C\in \cN^{(2)}$.
Then 
   there exists a unique
 $C'\in \cN^{(3)}$ and an element 
$g\in G^{(2)}$ such that   $C^g=C'$.
This implies that 
\begin{eqnarray*}
\cN_{torus}/G\cong\cN_{torus}^{(2)}/G^{(2)}\cong
\cN_{torus}^{(3)},\quad
\cN_{gen}/G\cong\cN_{gen}^{(2)}/G^{(2)}\cong
\cN_{gen}^{(3)}
\end{eqnarray*}
\end{lemma}
{\em Proof.}
Assume that  $C\in \cN_{torus}^{(1)}$,  $b+c=0$. Consider 
an element $g\in G_0^{(1)}$,
\[ g\inv=\left (
\begin{matrix}
1&0&0\\0&1&0\\ 1-a_3&0&a_3\end{matrix}\right)\]
The image 
$L_A^g$ is given by
$y-x+x a_3-a_3-bxa_3+ba_3=0$. Thus we can solve
the equation $a_3(1-b)-1=0$ in $a_3$ uniquely as
$a_3=1/(1-b)$ as  $b\ne 1$.  Thus $g\in G_0^{(1)}$ is unique if
it fixes the singular points pointwise and thus  $C'$ is also
unique.  It is easy to see that the stabilizer of
$\cN_{torus}^{(3)}$ is the cyclic group of order two 
generated by $\tau$, as $C'$ is even in $y$
(see the normal form below) and ${C'}^\tau=C'$ 
for any $C'\in \cN_{torus}^{(3)}$. Thus we have
$\cN_{torus}^{(2)}/G^{(2)}\cong
\cN_{torus}^{(3)}$. 

Consider the case $C\in \cN_{gen}^{(2)}$.
Then the images of the tangent cones at $A,B$ by the action of
$g$ are   given by
$y-x+x a_3-a_3-bxa_3+ba_3=0$ and 
$ y+x-xa_3+a_3-c x a_3+c a_3=0$ respectively.
 Assume that $b(C^g)=c(C^g)$. Then we need to have
$a_3(1-b)-1=a_3(-1-c)+1$,
which has a unique solution in $a_3$, if $(\star)$ $b-c-2\ne
0$.  Assume that $c^2+3c-bc+3-3b+b^2=0$ and $b-c-2=0$. 
Then we get $(b,c)=(1,-1)$ 
which is excluded as it corresponds to a non-reduced sextic. Thus 
the condition $(\star)$ is always satisfied.
Put $(b',c'):=(b(C^g),c(C^{g}))$. They  satisfy
 the equality 
 ${c'}^2+3{c'}-b'c'+3-3b'+{b'}^2=0$
and $b'=c'$. Thus we have  either $b'=c'=\sqrt{-3}$ or
$b'=c'=-\sqrt{-3}$. However in the second case, $(C^g)^\tau$
 belongs to  the first case. Thus 
$b'=c'=\sqrt{-3}$ and $C^g\in\cN_{gen}^{(3)}$ as desired. 
 \qed
\subsection{Normal forms of curves of a torus type}
In \cite{dual}, we have shown that a curve in $\cN_{torus}^{(1)}$
is defined by a polynomial 
$f(x,y)$ which
is expressed by  a sum $f_2(x,y)^3+s f_3(x,y)^2$
where $f_2(x,y)$ is a smooth conic passing through
$O,A,B$ and $f_3(x,y)=(y^2-x^2)(x-1)$.
\begin{proposition}
The direction of the 
 tangent cones at $O$, $A$ and $B$ are the same with the tangent
line of the conic
$f_2(x,y)=0$ at these points.
\end{proposition}
This is immediate as the multiplicity of $f_3(x,y)^2$ at
$O,A,B$ are 4. 
Assume that $C\in \cN_{torus}^{(3)}$, that is,
 the tangent cones of $C$ at $O$,  $A$ and $B$ are
given by 
$x=0$, $y-1=0$ and $y+1=0$ respectively. 
Thus   the
conic
$f_2(x,y)=0$ is also uniquely  determined as 
$f_2(x,y)=y^2+x^2-2x$. 
Therefore  $ \cN_{torus}^{(3)}$
 is one-dimensional and
it has the representation
\begin{eqnarray}\label{originalform}
C_s:\quad f_{torus}(x,y,s):=f_2(x,y)^3+sf_3(x,y)^2=0
\end{eqnarray}
For $s\ne 0,27,\infty$, $C_s$ is a sextic with three (3,4)-cusps, while
$C_{27}$ obtains a node. If
$g\in G^{(2)}$    fixes the tangent lines $y\pm 1=0$,
then
$g=e$ or $\tau$. As $C_s^\tau=C_s$, this implies that $C_s^g=C_s$. Thus 
$C_s\ne C_t$ if $s\ne t$.
\subsection{Normal form of sextics of a general type }
For the moduli
$\cN_{gen}$  of sextics of a general type, we start from the expression given
in
\S 4.1, \cite{dual}. 
We may assume $b=c=\sqrt{-3}$. Then the
parametrization is given by 
\begin{eqnarray*}
f_{gen}(x,y,s):=f_0(x,y)+s f_3(x,y)^2,\quad f_3(x,y)=(y^2-x^2)(x-1)
\end{eqnarray*}
where  $s$ is equal to 
$a_{06}$  in \cite{dual} and $f_0$ is 
the sextic given by 
\begin{eqnarray}\label{sqrt}
&f_0(x,y):= y^6+y^5(6\sqrt{-3}-6\sqrt{-3}x)+y^4
(35-76 x+38 x^2) 
\\&+y^3 (-24 \sqrt{-3} x+36 \sqrt{-3} x^2-12 
\sqrt{-3} x^3)  
 +y^2(-94 x^2+200 x^3-103 x^4)\notag\\&+
y (24 \sqrt{-3} x^3-42 \sqrt{-3} x^4+18 \sqrt{-3} x^5) 
+64 x^3-133 x^4+68 x^5\notag
\end{eqnarray}
Let $D_s:=\{f_{gen}(x,y,s)=0\}$ for each  $s\in \bfC$.
Observe that $D_0=\{f_0(x,y)=0\}$ is a sextic with three (3,4)-cusps
and of a general  type.
For the computation of dual curves  using  Maple V,  it is better to take the
substitution 
$y\mapsto y\sqrt{-3}$ to make the equation to be defined
over $\bfQ$.
Summarizing the discussion, we have 
\begin{theorem}\label{uniqueness}
The quotient  moduli space $\cN/G$ is one dimensional
and it has  two components.
\newline
\noindent (1)
 The component 
$\cN_{torus}/G$  has the normal forms 
  $C_s=\{f(x,y,s)=0\}$
 where 
$f(x,y,s)=f_2(x,y)^3+s f_3(x,y)^2$,
$f_2(x,y)=y^2+x^2-2x$ and   $f_3(x,y)=(y^2-x^2)(x-1)$.
The curve $C_{54}$ is
 a unique curve in $ \cN/G$
 which is self-dual. 
\newline
\noindent (2)
The  component $\cN_{gen}/G$ 
has the normal form:
$f_{gen}(x,y,s)=f_0(x,y)+sf_3(x,y)^2$ where $f_3$ is as above
and the sextic $f_0(x,y)=0$ is contained in $\cN_{gen}$. This component has
no self-dual curve.
\end{theorem}
{\em Proof of Theorem \ref{uniqueness}.}
We need only prove the assertion
for the dual curves. The proof will be done
 by a direct computation of dual curves using 
the method of \S 2, \cite{dual} and the above parametrizations.
We use Maple V for the practical computation.
Here is the recipe of
the proof. Let $X^*,Y^*,Z^*$ be the dual coordinates of $X,Y,Z$
and let $(x^*,y^*):=(X^*/Z^*,Y^*/Z^*)$ be the dual affine
coordinates.

(1) Compute the defining polynomials of the
dual  curves $C_s^*$ and $D_s^*$ respectively, using the method of
Lemma 2.4,
\cite{dual}. Put  $g_{torus}(x^*,y^*,s)$ 
and $g_{gen}(x^*,y^*,s)$ the defining polynomials.

(2) Let $G_\eps(X^*,Y^*,Z^*,s)$ be the homogenization of
$g_\eps(x^*,y^*,s)$, $\eps=torus$ or $gen$. Compute the discriminant
polynomials
$\De_{Y^*} G_\eps$ which is a homogeneous polynomial in $X^*,Z^*$ of
degree 30 (cf. Lemma 2.8,
\cite{flex-covering}). Recall that the multiplicity in $\De_{Y^*} G_\eps$
 of
the pencil $X^*-\eta Z^*=0$ passing through a singular point is 
generically given by $\mu+m-1$ where 
$\mu$ is the Milnor number and  $m$ is the multiplicity of the
singularity (\cite{dual}). Thus the contribution from  a (3,4)-cusp is  8.
Thus if $C_s^*$ has three (3,4)-cusps, it is necessary that 
$\De_{Y^*}(G)=0$ has three linear factors with multiplicity $\ge$
8.

(3-1)  For the  curves of a general type, an easy computation shows that it is
not possible to get a degeneration into a  sextic with 3 (3,4)-cusps
by the above reason.

(3-2) For the  curves of a torus type, we can see that $s=54$ is the only
 parameter such that $C_s^*\in \cN$.  Thus it is enough to show that
$C_{54}^*\cong C_{54}$.

(4) The dual curve $C_{54}^*$ of
$C_{54}$ is defined by
the homogeneous polynomial
\begin{eqnarray*}
G(X^*,Y^*,Z^*):=&128{X^*}^5 Z^*+1376{X^*}^4{Z^*}^2-192{X^*}^3 {Y*}^2
Z^* +4664 {X^*}^3 {Z^*}^3
-2{X^*}^2{Y^*}^4\\&
-1584{X^*}^2{Y^*}^2 {Z^*}^2+7090{X^*}^2{Z^*}^4+58
{X^*} {Y^*}^4 Z^*-3060 X^* {Y^*}^2 {Z^*}^3 \\&+5050 X^*
{Z^*}^5+{Y^*}^6+349 {Y^*}^4 {Z^*}^2-1725 {Y^*}^2 {Z^*}^4+1375
{Z^*}^6
\end{eqnarray*}
We can see that ${C_{54}}^*$ is
isomorphic to $C_{54}$ as $(C_{54}^*)^A=C_{54}$  where 
\[ A=\begin{pmatrix}
-4/3&0&-5/3\\0&1&0\\-5/3&0&-13/3\end{pmatrix}
\]

\subsection{Involution $\tau$ on $C_{54}$}
For a later purpose, we change the coordinates of $\bfP^2$
 so that the three cusps  of $C_s$ are at 
$O_Z:=(0,0,1),O_Y:=(0,1,0),O_X:=(1,0,0)$.
A new normal form
in the  affine space is given by
$ C_s:~f_2(x,y)^3+sf_3(x,y)^2=0$ where
$f_2(x,y):=xy-x+y$ and $f_3(x,y):=-xy$.
In particular,
$C_{54}$ is defined by 
\begin{eqnarray}\label{sextic-eq}
f(x,y)=(xy-x+y)^3+ 54x^2y^2=0
\end{eqnarray}
In this coordinate, $C_{54}^* $ is defined by 
\begin{eqnarray*}
&-28 y^3-17 x^4 y^2-17 x^2 y^4-28 x^3 y^3-2 y^5+1788 x^3 y 
+1788 x^2 y-17 y^4-17 x^4\\&
+262 x y+1788 x^2 y^3-1788 x y^2 
-262 x y^4+1788 x y^3-1788 x^3 y^2-8166 x^2 y^2+28 x^3\\&+262 x^4 y 
-2 x^5 y-2 x y^5+1-17 y^2-17 x^2+2 x^5+2 x-2 y+x^6+y^6=0
\end{eqnarray*}
It is easy to see that
$(C_{54}^*)^{A_1}=C_{54}$ where 
\[
A_1:=\begin{pmatrix}
-1/3& 7/3& -1/3\\ 7/3& -1/3& 1/3\\ -1/3&1/3&-7/3\end{pmatrix}
\]
For a given $A\in \GL(3,\bfC)$, we denote the automorphism defined by
 the right multiplication of $A$ by $\vphi_A$.
Let $F(X,Y,Z)$ be the homogenization of $f(x,y)$.
Then the Gauss map  $\dual_{C_{54}}:C_{54}\to C_{54}^*$
 is defined by 
\[
\dual_{C_{54}}(X,Y,Z)= (F_X(X,Y,Z),F_Y(X,Y,Z),F_Z(X,Y,Z))\]
 where
$F_X,F_Y,F_Z$ are partial derivatives.
We define an isomorphism $\tau:~C_{54}\to
C_{54}$ by the composition 
 $\vphi_{A_1}\circ \dual_{C_{54}} $.
Then  $\tau$ is the
restriction of the rational mapping:
$\Psi:\bfC^2 \to \bfC^2$, $(x,y)\mapsto (x_d,y_d)$ and 
\begin{eqnarray}\label{tau-sextic}\begin{cases}
&x_d :=\frac{
-y^3+4x^2-x^2 y^3+4 x^3 y^2-8 x^3 y-4 x^2 y^2-8 x y-4 x y^2-2 x  y^3
+109 x^2 y+4 y^2+4 x^3}{ -4 y^3+x^2-4 x^2 y^3+4 x^3 y^2-8 x^3 y
-109 x^2  y^2-2 x y-4 x  y^2-8 x y^3+4 x^2 y+y^2+4 x^3}\\
&y_d :=-\frac{
-4 y^3+4 x^2-4 x^2 y^3+x^3 y^2-2 x^3 y-4 x^2 y^2-8 x  y-109 x y^2
-8 x y^3+4 x^2 y+4 y^2+x^3}{-4 y^3+x^2-4 x^2 y^3+4 x^3 y^2-8 x^3 y
-109 x^2 y^2-2 x y-4 x y^2-8 x y^3+4 x^2  y+y^2+4  x^3}
\end{cases}\end{eqnarray}
Observe that $\tau$ is defined over $\bfQ$.
$C_{54}$ has three flexes of order 2 at 
$F_1:=(1,-1/4,1), ~F_2:=(1/4,-1,1),~F_3:=(4,-4,1)$ and 
$\tau$ exchanges flexes and cusps:
\begin{eqnarray}
\begin{cases}&\tau(O_X)=F_1,  \tau(O_Y)=F_2,\tau(O_Z)=F_3,\\
&\quad  \tau(F_1)=O_X,  \tau(F_2)=O_Y,  \tau(F_3)=O_Z\end{cases}
\end{eqnarray}
Furthermore we assert that 
\begin{proposition}\label{involution6}
The morphism $\tau$ is an involution on $C_{54}$.
\end{proposition}
{\em Proof.}
By the definition of $\tau$ and Lemma \ref{GaussMap} below,
we have ($C:=C_{54}$):
\[\tau\circ \tau=(\vphi_{{}^tA_1\inv}\circ\dual_C)^2
=(\dual_{C^{A_1}}\circ\vphi_{A_1})\circ(\vphi_{{}^tA_1\inv}\circ\dual_C)
=\id\]
as $A_1$ is a symmetric matrix.\qed

Let $C$ be a  given irreducible curve in $\bfP^2$ defined by a 
homogeneous polynomial $F(X,Y,Z)$
and let 
$B\in \GL(3,\bfC)$. Then $C^B$ is defined by
$G(X,Y,Z):=F((X,Y,Z)B\inv)$.
\begin{lemma}\label{GaussMap} Two curves $(C^B)^*$ and $(C^*)^{{}^tB\inv}$
coincide and the following
diagram commutes.
\[
\begin{matrix}
C&\mapright{\dual_C}& C^*\\
\mapdown{\vphi_{B}}&&\mapdown{\vphi_{{{}^tB}\inv}}\\
C^B&\mapright{\dual_{C^B}}&(C^B)^*
\end{matrix}
\]
\end{lemma}
{\em Proof.}
The first assertion is the same as Lemma 2, \cite{dual}.
The second  assertion follows from the following equalities. 
Let $(a,b,c)\in
C$.
\begin{eqnarray*}
&\dual_{C^B}(\vphi_B(a,b,c))=
(G_X(\vphi_B(a,b,c)),G_Y(\vphi_B(a,b,c)),G_Z(\vphi_B(a,b,c)))\\
&\quad =
(F_X(a,b,c),F_Y(a,b,c),F_y(a,b,c)){}^tB\inv=
\vphi_{{}^tB\inv}(\dual_C(a,b,c)\qed
\end{eqnarray*}

 In  section 5, we will show that $\tau$ is expressed in a simple
form as a
cubic curve.

\vspace{.3cm}
\section{Structure of elliptic fibrations}

We consider the elliptic fibrations corresponding to the above 
normal forms. For this purpose, we first take a linear change  of
coordinates so that three lines defined by $f_3(x,y)=0$ changes
into 
  lines
$X=0$, $Y=0$ and $Z=0$.
The corresponding three cusps are now at 
$O_Z=(0,0,1),O_Y=(0,1,0), O_X=(1,0,0)$ 
in $\bfP^2$. Then we take the quadratic 
transformation which 
 is a birational mapping of $\bfP^2$ defined by 
$(X,Y,Z)\mapsto (YZ,ZX,XY)$.
Geometrically this is the composition of   blowing-ups at
$O_X,O_Y,O_Z$ and  then the blowing down of
three lines which are
strict transform of 
$X,Y,Z=0$.
It is easy to see that our sextics are transformed 
into smooth cubics for which 
$X=0$, $Y=0$ and $Z=0$  are tangent lines of the flex points.
Those flexes are the image of the (3,4)-cusps.
We take a 
  linear change of coordinates so that the flex on $Z=0$ 
is moved at $O:=(0,1,0)$ with the tangent $Z=0$. Then
 the corresponding families are described by the  families given by
$\{h_{torus}(x,y,s)=0;s\in \bfP^1\}$ and 
$\{h_{gen}(x,y,s)=0,s\in \bfP^1\}$ where 
\begin{eqnarray*}
\begin{cases}
&C_s:\quad h_{torus}(x,y,s):= x^3-\frac1 4 s (x-1)^2+s y^2,\quad\qquad\qquad\\
&D_s:\quad h_{gen}(x,y,s) :=
-8 x^3+1+(s+35)y^2-6 x^2+3x\\&
\qquad \qquad -6\sqrt{-3}y-3
\sqrt{-3}x -6\sqrt{-3}x^2-12\sqrt{-3} xy+(s-35)xy
\end{cases}
\end{eqnarray*}
Let  $H_{\eps}(X,Y,Z,S,T)=h_\eps(X/Z,Y/Z,S/T)Z^3T$ for  $\eps=torus,~gen$.
  We consider the elliptic surface associated to
 the canonical projection
 $\pi:S_{\eps} \to \bfP^1$
where $S_{\eps}$ is the hypersurface in 
$\bfP^1\times\bfP^2$ which is defined by
$H_\eps(X,Y,Z,S,T)=0$.

\vspace{.2cm}
 Case I. Structure of   $S_{torus}\to \bfP^1$.
 For simplicity, we use
the affine coordinate $s=S/T$ of
$\{T\ne 0\}\subset
\bfP^1$ and denote $\pi\inv(s)$ by $C_s$.
 We see that the singular fibers are $
s=0,27,\infty$. $C_\infty$ consists of  three lines,  isomorphic
to $I_3$  in  Kodaira's notation, \cite{Kodaira}.
$C_{27}$ obtains a node and this fiber is  denoted by $I_1$  in
\cite {Kodaira}. The fiber $C_0$ is a line with multiplicity
3. The surface $S_{torus}$ has three singular points on the fiber $C_0$
at $(X,Y,Z)=(0,1/2,1),(0,-1/2,1),(0,1,0)$. Each singularity is
an $A_2$-singularity. We
take minimal resolutions at these points. At each point, we
need two $\bfP^1$ as exceptional divisors and let 
$p:{\wtl {S}}_{torus}\to S_{torus}$ be the resolution map.
The composition 
$\wtl\pi:=\pi\circ p:{\wtl{S}}_{torus}\to\bfP^1$  is the
corresponding elliptic surface. Now it is easy to see that
$\wtl{ C_0}:={\wtl \pi}\inv (0)$ is a singular fiber with 7
irreducible components, which is denoted by $IV^*$ in
\cite{Kodaira}. Here we used the  following lemma. 
The elliptic surface
$\wtl{S}_{torus}$ is rational and denoted by $X_{431}$ in
\cite{Miranda}.

Assume that the surface 
$V:=\{(s,x,y)\in \bfC^3;f(s,x,y)=0\}$ has an $A_2$ singularity at the origin
where  $f(s,x,y):=sx+y^3+sx \cdot h(s,x,y)$ where $h(O)=0$. 
Consider the minimal resolution
$\pi:\wtl V \to V$ and let 
$\pi\inv(O)=E_1\cup E_2$.
It is well-known  that $E_1\cdot E_2=1$ and $E_i^2=-2$
for $i=1,2$. 
\begin{lemma}\label{A2-data}Consider a linear form $\ell(s,x,y)=as+bx+cy$
and let $L'$  be the strict transform of $\ell=0$ to $\wtl V$.

\noindent
(1) Assume that $b=c=0$ and $a\ne 0$. Then 
$(\pi^*\ell)=3 L'+2 E_1+E_2$ and $L'\cdot E_1=1$ and 
$L'$ does not intersect with $E_2$, under a suitable ordering
of $E_1$ and $E_2$.

\noindent
(2) Assume that $abc \ne 0$. Then we have 
$(\pi^*\ell)= L'+ E_1+E_2$ and $L'\cdot E_i=1$ for $i=1,2$.
\end{lemma}
The proof is immediate from a direct computation.

\vspace{.3cm}
Case II.  Structure of  $S_{gen}\to \bfP^1$. 
Now consider the
elliptic surface
$S_{gen}$.   Put 
$D_s=\pi\inv(s)$.
The singular fibers are at $s=-35, -53+6\sqrt{-3}$,
$-53-6\sqrt{-3}$ and $s=\infty$. The fiber $s=\infty$
is already $I_3$ and $S_{gen}$ is smooth on this fiber.
On the other hand,  $S_{gen}$
has a $A_2$-singularity on each  fiber $D_s$,  $s=-35, -53+6\sqrt{-3}$,
$-53-6\sqrt{-3}$.
Let $p: \wtl{S}_{gen}\to S_{gen}$ be the the minimal
resolution map and  we consider the composition
  $\wtl\pi:=\pi\circ p:\wtl{S}_{gen}\to \bfP^1$  as above. Then
using (2) of  Lemma
\ref{A2-data},
 we see that 
$\wtl \pi: {\wtl {S}}_{gen}\to\bfP^1$ 
has four singular fibers and each of them is $I_3$ in the notation 
\cite{Kodaira}. This elliptic surface is also rational and denoted as 
$X_{3333}$ in \cite{Miranda}.

\vspace{.3cm}
\section{Torsion group of $C_s$ and $D_s$}
Consider an elliptic curve  $C$ defined
over a  number field $K$ by a Weierstrass short normal
form:
$ y^2=h(x),\quad h(x)= x^3+Ax+B$. 
The j-invariant is defined by $j(C)=-1728(4A)^3/\De$ with
 $\De=-16(4A^3+27B^2)$. We study  the torsion group of the
Mordell-Weil group of 
$C$ which we denote by $C_{tor}(K)$ hereafter. 

Recall that a point of order 3 is geometrically
a flex point of the complex curve $C$ (\cite{Silverman})
 and its locus is defined by
$\cF(f):=f_{x,x}f_y^2-2f_{x,y}f_xf_y+f_{y,y}f_x^2=0$ where $f(x,y)$ is 
the defining polynomial of $C$
(\cite{flex-covering}). In our case, 
$\cF(f)= 24xy^2-18x^4-12x^2 A-2A$.
 The unit of the group is given by the point
at infinity $O:=(0,1,0)$ and the inverse of
$P=(\al,\be)\in C$ is given by 
$(\al,-\be)$ and we denote it by $-P$.
For a later purpose, we prepare two easy propositions.
Consider a line $L(P,m)$ passing through $-P$ defined by
$y=m(x-\al)-\be$. The $x$-coordinates of two other intersections
with $C$ are the solution of 
$ q(x):=f(x,m(x-\al)-\be)/(x-\al)$
which is a polynomial of degree 2 in $x$. Let $\De_x q$ be the 
discriminant of $q$ in $x$. Note that  $\De_x q$
is a polynomial in $m$.

\vspace{.2cm}\noindent
{ (A) When does  a point $Q\in C$ exist such that $2Q=P$.}
\nl
Assume that a $K$ point $Q=(x_1,y_1)$ satisfies $2Q=P$. Geometrically 
this implies that the tangent line $T_Q C$ passes through $-P$.
\begin{proposition}\label{2Q} There exists 
a $K$-point $Q$ with $2Q=P$
 if and only if
$m$ is a $K$-solution of $\De_x q (m)=0$ and 
$x_1$ is the multiple solution of $q(x)=0$.
If $P$ is a flex point, $\De_x q(m)=0$ contains a canonical
solution which corresponds to the tangent line at $P$
and $m=-f_x(\al,\be)/f_y(\al,\be)$. For any $K$-solution $m$
with $m\ne -f_x(\al,\be)/f_y(\al,\be)$, 
the order of  $Q$ is equal to  $2\cdot \order~ P$.
\end{proposition}

\vspace{.2cm}\noindent
{(B) When does   a point $Q\in C$ exist such that $3Q=P$.}\nl
Assume that a $K$-point $Q=(x_1,y_1)$ satisfies $3Q=P$. 
Put $Q':=2Q$ and put $Q'=(x_2,y_2)$. Let $T_Q C$ be the
tangent line at $Q$. Then $T_Q C$ intersects  $C$ at $-Q'$.
Then $-3Q$ is the third intersection of  $C$ and the line $L$ 
which passes through $Q,Q'$.  Thus three points $-P,Q,Q'$ are 
colinear. Write $L$ as $y=m(x-\al)-\be$.
Then $x_1,x_2$ are the solutions of $q(x)=0$.
Thus we have 
\begin{eqnarray}
x_2=-\coeff(q,x)/\coeff(q,x^2)-x_1,\quad y_1=m(x_1-\al)-\be
\end{eqnarray}
where 
$\coeff(q,x^i)$ is the coefficient of $x^i$ in $q(x)$.
Let $L_Q(x,y)$ be the linear form defining $T_Q C$ and 
let $R(x)$ be  the resultant  of $f(x,y)$ and $L_Q(x,y)$
in $y$. 
 Put
$R_1(x):=R(-\coeff(q,x)/\coeff(q,x^2)-x)$.
Then
by the above consideration,  $x=x_1$ is a common solution of
 $q(x)=R_1(x)=0$.
Let $R_2(m)$ be the resultant of $q(x)$ and $R_1(x)$.
Note that if  $\De_x q(m)=0$, $L$ is tangent to $C$ at $Q$ and $R_2(m)=0$.
In this case, $2Q=P$.
\begin{proposition}\label{3Q}
Assume that there exists a $K$-point $Q$ with $3Q=P$ and 
$\order~ Q=3\cdot\order~ P$ and let $m$ be as above.
Then  
$R_2(m)=0$  and  $\De_x q(m)\ne 0$.
Moreover  $x_1$ is given as a common solution of $q(x)=R_1(x)=0$.
\end{proposition}
Actually one can show that $R_2(m)$ is divisible by $(\De_x q)^2$. 
 \subsection{  Cubic family associated with sextics of a
torus type}
 We have observed that the family 
$C_s$ for $s\in \bfQ$ is defined over
$\bf Q$. 
 First, recall that $C_s$ is defined by 
\begin{eqnarray}
C_s: x^3-\frac 1 4 s (x-1)^2+s y^2=0
\end{eqnarray}
 and  
the Weierstrass normal form is given by
$C_s: ~y^2=x^3+a(s)x+b(s)$ where 
\begin{eqnarray}\label{coeff-torus}
 a(s)=-\frac 1{48}s^4+\frac 12 s^3,~b(s)=-\frac{1}{24}s^5+\frac
14 s^4 +\frac 1{864} s^6\end{eqnarray}
Put $\Si:=\{0,27,\infty\}$. This corresponds to singular fibers.
We have two sections of  order 3:
$s\mapsto (\frac 1{12} s^2, \pm \frac 12 s^2)$.
Put $P_1:=(\frac 1{12} s^2, \frac 12 s^2)$.
Thus the torsion group is at least $\bfZ/3\bfZ$. 
By
\cite{Mazur}, 
 the possible torsion group which has an element of order 3 
is one of  $\bfZ/3\bfZ,~\bfZ/6\bfZ$,
 $\bfZ/2\bfZ+\bfZ/6\bfZ$, $\bfZ/9\bfZ$ or $\bfZ/12\bfZ$.
 The j-invariant of $C_s$ is given by 
\begin{eqnarray}
j(C_s):=j_{torus}(s),\quad j_{torus}(s):=s(s-24)^3/(s-27)
\end{eqnarray}

\vspace{.3cm}\noindent
(1) Assume that $(C_s)_{tor}(\bfQ)$ has an
element of order $6$, say $P_2:=(\al_2,\be_2)\in C_s\cap \bfQ^2$. We may assume
that
$P_2+P_2=P_1$. By Proposition \ref{2Q}, this implies that  
 $x=\al_2$ must be the multiple solution of 
$$q(x):=s^4-36 s^3-72 m s^2-6 x s^2-6 s^2 m^2+72 m^2 x-72 x^2=0$$
As $-f_x(-P_1)/f_y(-P_1)=-s/2$, we must  have  $m\ne -s/2$ and thus 
\begin{eqnarray}\label{order6-eq}
\De_x'q:=\De_x q/(2m+s)=s^3-32 s^2-2 m s^2-4 m^2 s+8 m^3=0\end{eqnarray}
The curve $\De_x'(q)=0$ is a rational curve 
and we can parametrize  $\De_x' q =0$
 as $s=\vphi_6(u)$, $m=\vphi_6(u)u$ where 
\begin{eqnarray}
\vphi_6(u):=32/(1+2u)(2u-1)^2\end{eqnarray}
The point $P_2$ is parametrized as
\begin{eqnarray}
P_2=( \frac{128}3 \frac{-1+12u^2}{(2 u+1)^2 (-1+2 u)^4},
\frac{ 512(6 u+1)}{(-1+2u)^5 (2u+1)^2})
\end{eqnarray}
where $u\in \bfQ$.
We put
$A_6:=\{s=\vphi_6(u);u\in \bfQ\}$ and
$\Si_6:=\vphi\inv(\Si)$.
Note that $\Si_6=\{-1/2,1/2,5/6,-1/6\}$.

\vspace{.3cm}\noindent
(1-2)
Assume that we are given $s=\vphi(u)$ and we consider the case when 
(\ref{order6-eq}) has three rational solutions in $m$ for a fixed $s$.
This  is the case if $\vphi_6(u)=\vphi_6(v)$
has two rational solutions different from $u$. This is also
equivalent to
$(C_s)_{tor}(\bfQ)$ has $\bfZ/2\bfZ+\bfZ/2\bfZ$ as a subgroup.
The equation  is given by the conic
\begin{eqnarray}
Q:\quad 4 u^2-2 u+4 u v-1-2 v+4 v^2=0\end{eqnarray}
By an easy computation, $Q$  is rational and it has a parametrization as
follows.
\begin{eqnarray}
u=\vphi_2(r):=\frac {-36+5 r^2}{6(12+r^2)},\quad
v(r):=-\frac 16 \frac{(r^2+24 r-36)}{(12+r^2)}
\end{eqnarray}
The generators are $P_2$ of order 6 and $R=(\ga,0)$   of order 2
where 
\[\ga:=-\frac{81}{4}\frac{(r^4-48r^3+72r^2-432)(12+r^2)^4}
{(r^2-36)^4r^4}\] 
Put $\vphi_{6,2}(r):=\vphi_6(\vphi_2(r))$, which is given explicitly as
\[\vphi_{6,2}(r)=27(12+r^2)/{r^2(r-6)^2(r+6)^2}\]
We define a subset $A_{6,2}$ of $A_6$ by the image $\vphi_{6,2}(\bfQ)$.
Put $\Si_{6,2}:=\vphi_{6,2}\inv(\Si)$. It is given by
 $\Si_{6,2}=\{0,\pm 2,\pm 6\}$.

\vspace{.3cm}\noindent
(2) Assume that there exists  a rational point $P_3=(\al_3,\be_3)$
of order 9
such that $3P_3=P$. By Proposition \ref{3Q}, this is the case if and
only if 
\begin{eqnarray*}
&R_3(m,s):=512 m^9+768 m^8 s-512 m^6
s^3-1536 m^6 s^2-192 s^4 m^5 \\& -6144 m^5 s^3-6528 m^4 s^4+96 s^5
m^4-12288 m^3 s^4 -2048 m^3 s^5+64 s^6 m^3 
+480 s^6 m^2\\&-15360 s^5 m^2 -6144 s^6 m 
+384 s^7 m-6 s^8 m+56 s^8-512 s^6-768 s^7-s^9 =0
\end{eqnarray*}
has a rational solution. Here $R_3$ is $R_2/(\De_x q)^2 (s+2m) s^4$
up to a constant multiplication.
Again we find that the curve
 $\{(m,s)\in \bfC^2; R_3(m,s)=0\}$ is rational and we can parametrize
this curve as $s=\vphi_9(t),~m=\psi_9(t)$ where
\begin{eqnarray}\begin{cases}
&\vphi_9(t):= -\frac 1 8\frac{ (-1+9 t^2-3 t+3 t^3)^3}
{t^3 (t-1)^3 (t+1)^3}\\
&\psi_9(t):=\frac 1{16} \frac{(-1+9 t^2-3 t+3 t^3)^2 (-t+t^3+1+7
t^2)}{t^3 (t-1)^3 (t+1)^3}\end{cases}
\end{eqnarray}
The generator $P_3=(\al_3,\be_3)$ is given by
\[\begin{cases}
&\al_3=\frac{1}{768}\frac{
(1-18 t+15  t^2-12 t^3+15 t^4+30 t^5+33 t^6) (9 t^2-1+3 t^3-3 t)^4}
{(t-1)^6 (t+1)^6 t^6}\\
&\be_3=-\frac{1}{512}\frac{(1+3 t^2) (9 t^2-1+3 t^3-3 t)^6}
{(t-1)^5 (t+1)^7 t^8}\end{cases}
\]
We put 
$ A_9:=\{\vphi_9(t);t\in \bfQ\}$ and 
$\Si_{9}:=\vphi_9\inv(\Si)=\{0,1,-1\}$.

\vspace{.3cm}\noindent
(3)
 Assume that $ s\in A_6$  and $(C_s)_{tor}(\bfQ)$ has an element
$P_4=(\al_4,\be_4)\in C_s\cap \bfQ^2$ of order $12$. 
Then we may assume that
$P_4+P_4=P_2$. This implies that the tangent line at $P_4$ passes
through
$-P_2$. Write this line as $y=n(x-\al_2)-\be_2$.
By the same discussion as above,  the
equality $\Ga(n_1,u)=0$ holds where $\Ga$ is 
the polynomial defined by
\begin{eqnarray}
&\Ga(u,n_1):=-786432 u^4-98304 n_1 u^3-524288 u^3+393216
u^2-16384 n_1 u^2\\& -3072 n_1^2 u^2+131072 u+24576 n_1 u+4096
n_1+16384+256 n_1^2+n_1^4\notag
\end{eqnarray}
 and  $n=n_1/(2u+1)(2u-1)^2$.  Again we find that 
$\Ga=0$ is a rational curve and 
 we have a parametrization:
$u=u(\nu)$ and $n_1=n_1(\nu)$ where
\begin{eqnarray}
&u(\nu)=-\frac 12\frac{ (\nu^4+2 \nu^2+5)}{(\nu^4-6 \nu^2-3)},
\quad
n_1(\nu)=-16 \frac{(2 \nu^2-4 \nu^3-4 \nu+\nu^4-3)}{(\nu^4-6 \nu^2-3)}\\
&s=\vphi_{12}(\nu):=\vphi_6(u(\nu)),\quad  \vphi_{12}(\nu):=
  -\frac{(\nu^4-3-6 \nu^2)^3}{(\nu-1)^4(1+\nu)^4(1+\nu^2)}
\end{eqnarray}
The generator of the torsion group $\bfZ/12\bfZ$ is 
$P_4=(\al_4,\be_4)$ where
 \begin{eqnarray*}\begin{cases}
&\al_4:=
\frac {1}{12}\frac{(\nu^8-12 \nu^7+24 \nu^6-36 \nu^5+42 \nu^4+12 \nu^3+36 \nu-3)
(\nu^4-6 \nu^2-3)^4}{(\nu-1)^8 (\nu+1)^8 (\nu^2+1)^2}\\
&\be_4:=
-\frac 12\frac{(\nu^4-6 \nu^2-3)^6 \nu (\nu^2+3)}{(\nu-1)^7 (\nu+1)^{11} (\nu^2+1)^2}
\end{cases}\end{eqnarray*}
 We put
$A_{12}:=\{\vphi_{12}(\nu);\nu\in \bfQ\}$. By definition,
$A_{12}\subset A_6$. The singular fibers $\Si_{12}:=\vphi\inv(\Si)$ is given
by $\{0,\pm 1\}$.
Summarizing the above discussion, we get
\begin{theorem}\label{torus-torsion} 
The j-invariant is given by $j_{torus}(s)=s(s-24)^3/(s-27)$ and 
the Mordell-Weil torsion group   of $C_s$ is given as follows.
\[(C_s)_{tor}(\bfQ)=\begin{cases}
&\bfZ/3\bfZ,\quad  s\in \bfQ-A_6\cup A_9\cup \Si\\
&\bfZ/6\bfZ,\quad s=\vphi_6(u)\in A_6-A_{6,2}\cup A_{12},~u\in \bfQ-\Si_6\\
&\bfZ/6\bfZ+\bfZ/2\bfZ,\quad  s=\vphi_{6,2}(r)\in A_{6,2},~r\in
\bfQ-\Si_{6,2}\\ 
&\bfZ/9\bfZ,\quad s=\vphi_9(t)\in A_9,~t\in \bfQ-\Si_{9}\\
& \bfZ/12\bfZ,\quad s=\vphi_{12}(\nu)\in A_{12},~\nu\in
\bfQ-\Si_{12}
\end{cases}\]
\end{theorem}
\subsection{Comparison with Kubert family}
In \cite{Kubert}, Kubert gave parametrizations of the moduli of elliptic curves
defined over $\bfQ$ with given torsion groups which have an element of
order $\ge 4$. His family starts with the normal form:
\begin{eqnarray}
E(b,c): y^2+(1-c)xy-by=x^3-bx^2
\end{eqnarray}
We first eliminate  the linear term of $y$ and then the
coefficient of $x^2$.  Let $K_w (b,c)$ be
 the Weierstrass short normal form,  which is obtained in this
way. The j-invariant is given by
\[
j(E(b,c))=\frac {(1-8 b c^2-8 c b-4 c+16 b+6 c^2+16 b^2-4 c^3+c^4)^3}
{b^3 (3 c^2-c-3 c^3-8 b c^2+b-20 c b+c^4+16 b^2)}\]
For a given elliptic curve $E$ defined over $K$
with Weierstrass normal form
$ E: y^2=x^3+ax+b$ and a given $k\in K$, the change of coordinates
$x\mapsto x/k^2,y\mapsto y/k^3$ changes the normal form into 
$y^2=x^3+ak^4 x+ bk^6$. We denote this operation by 
$\Psi_k(E)$.

1. Elliptic curves with the torsion group $\bfZ/6\bfZ$.
This family is given by a parameter $c$ with $b=c+c^2$.
   
2. Elliptic curves with the torsion group $\bfZ/6\bfZ+\bfZ/2\bfZ$.
This family is given by a parameter $c_1$ with $b=c+c^2$
and $c=(10-2c_1)/(c_1^2-9)$.

3. Elliptic curves with the torsion group $\bfZ/9\bfZ$.
The corresponding parameter is $f$ and 
$b=cd,c=fd-f, d=f(f-1)+1$.

4. Elliptic curves with the torsion group $\bfZ/12\bfZ$.
The corresponding parameter is $\tau$ and 
$b=cd,c=fd-f, d=m+\tau, f=m/(1-\tau)$
and $m=(3\tau-3\tau^2-1)/(\tau-1)$.

\begin{proposition}
Our family $C_{\vphi_6(u)},
C_{\vphi_{6,2}(r)},C_{\vphi_9(t)},C_{\vphi_{12}(\nu)}$ are
equivalent to the respective Kubert families. More explicitly,  we 
take the following change of parameters to make   their
j-invariants  coincide  with those of Kubert and then we take the 
change of coordinates of type  $\Psi_k$  to make the Weierstrass
short normal forms to be identical with $K_w(x,y)$.

\begin{enumerate}

\item For $C_{\vphi_6(u)}$, take
$u=-(c-1)/2(3c+1)$ and 
$ k=c^2(c+1)/(3c+1)^2$.
\item For $C_{\vphi_{6,2}(r)}$, take
$r= -12/(c_1-3)$ and 
$k=4 (-5+c_1)^2 (c_1-1)^2/(c_1^2-6 c_1+21)^2$
\newline
$/ (c_1-3) (c_1+3)$.
\item For $C_{\vphi_9(t)}$, take  $t=-f/(f-2)$ and 
$ k=f^3(f-1)^3/(f^3-3f^2+1)^2$.
\item For $C_{\vphi_{12}(\nu)}$, take
$\nu=-1/(2\tau-1)$ and 
$k=(\tau-1) \tau^4 (-2 \tau+2 \tau^2+1)
 (-1+2 \tau)^2/$\nl
$(6 \tau^4-12 \tau^3+12 \tau^2-6 \tau+1)^2$.
\end{enumerate}
\end{proposition}
We omit the proof as the assertion  is immediate from  a direct computation.

\subsection{Involution on $C_{54}$}\label{s=54}
We consider again  the self dual curve $C:=C_{54}$ (see \S 3).
The Weierstrass normal form
is $y^2=x^3-98415x+11691702$. 
Note that $54\in A_6-A_{12}\cup A_{6,2}\cup
\Si$. In fact, $54=\vphi_6(1/6)$  and $54\notin
A_{12}\cup A_{6,2}$. 
The j-invariant is 54000 and the torsion group 
$C_{tor}(\bfQ)$ is $\bfZ/6\bfZ$ and 
the generator is given by $P=(-81,4374)$. Other rational points are
$2P=(243,-1458),3P=(162,0),4P=(243,1458),5P=(-81,-4374)$, and 
$O=(0,1,0)$ (= the point at infinity). Recall that 
$C$ has an involution $\tau$ which is defined by (\ref{tau-sextic})
 in \S 3.
To distinguish our original sextic and cubic, we put 
\[
C^{(6)}: (xy-x+y)^3+54x^2y^2=0,\quad C^{(3)}:y^2=x^3-98415x+11691702\]
The identification $\Phi:C^{(3)}\to C^{(6)}$ is given by the rational mapping:
\[\Phi(x,y)=(-2916/(27 x-5103-y),2916/(y+27x-5103))\]
and the involution $\tau^{(3)}$ on $C^{(3)}$ is given by the
composition $\Phi\inv\circ\tau\circ \Phi$.
After a boring computation,  $\tau^{(3)}$ is reduced to an extremely simple form
in the Weierstrass normal form  and it
is given
by
$\tau^{(3)}(x,y)=(p(x,y),q(x,y))$ where 
\begin{eqnarray}\label{tau-cubic}
 p(x,y) :=81\frac{2x-567}{x-162}\quad
q(x,y) :=-19683\frac{y}{(x-162)^2}
\end{eqnarray}
Note that $C$ has another  canonical
involution $\iota$ which is an automorphism defined by
$\iota: (x,y)\mapsto (x,-y)$. We can easily check that 
$\tau^{(3)}\circ\iota=\iota\circ \tau^{(3)}$.  Note
that 
$\tau^{(3)}(P)=2P,\tau^{(3)}(2P)=P,\tau^{(3)}(3P)=O,\tau^{(3)}(O)=3P,\tau^{(3)}(4P)=5P,\tau^{(3)}(5P)=4P$. Let
$\eta:C\to C$ be the translation by the 2-torsion
element
 $3P$ i.e., $\eta(x,y)=(x,y)+(162,0)$.
It is easy to see that  $\tau^{(3)}$ is the
 composition $\iota\circ \eta$. That 
is $\tau^{(3)}(x,y)=(x,-y)+(162,0)$ where the addition is the addition
by the  group structure of $C_{54}$. Thus 
\begin{theorem}The involution $\tau$ on sextics $C^{(6)}$
is equal to
the involution $\tau^{(3)}$ on $C^{(3)}$
which is defined by (\ref{tau-cubic}) and 
it is also equal to 
$(x,y)\mapsto (x,-y)+(162,0)$.
\end{theorem}

\vspace{.5cm}
\subsection{ Cubic family associated with  sextics of a general
type} We consider the family of elliptic $D_s$ curves associated to
the moduli of sextics of a general type with three (3,4)-cusps.
Recall that $D_s$ is defined by the equation:
\begin{eqnarray*}
&D_s:\quad -8 x^3+1+sy^2+35 y^2-6 x^2+3x-6\sqrt{-3}y-3 \sqrt{-3}x
\\&\qquad -6 \sqrt{-3}x^2-12\sqrt{-3} xy+(s-35)xy=0
\end{eqnarray*}

This family is defined over $\bfQ(\sqrt{-3})$. 
 We change this polynomial
into a Weierstrass normal form  by the usual process killing 
the coefficient of
$y$ and then by killing the coefficient of $x^2$.
 A  Weierstrass normal forms is given by $y^2=x^3+a(s)x+b(s)$
where  
\begin{eqnarray}\label{coeff-nontorus}\ \begin{cases}
&a(s):=-\frac{1}{768} (s+47) (s+71) (s^2+70 s+1657)\\
&b (s):=
\frac{1}{55296} (s^2+70 s+793) (s^4+212 s^3+17502 s^2\\
&\qquad\qquad \qquad \qquad \qquad\qquad +648644 s+9038089)
\end{cases}\end{eqnarray}
The singular fibers are $s=-35, -53+6\sqrt{-3}$,
$-53-6\sqrt{-3}$ and $s=\infty$. Put $\Si=\{-35,-53+\pm 6 \sqrt{-3},\infty\}$.
In this
section, we consider the Modell-Weil torsion over the quadratic number field
$\bfQ(\sqrt{-3})$.
 First we  observe
that this family has  8  sections of order three
$\pm P_{3,i},i=1,\dots,4$ where 
 $P_{3,i}$ are given by 
\begin{eqnarray}
&P_{3,1}:=(x_{3,1},y_{3,1} ),~\begin{cases} &x_{3,1 }:=
5041/48+71 s/24  +  s^2/48 \\ &y_{3,1}:=2917/4+53s/2
 +s^2/4 \end{cases} \\ 
 &P_{3,2}:=(x_{3,2} ,y_{3,2} )\begin{cases}&x_{3,2 }:=
-2209/16-47 s/8  -s^2/16  \\ & y_{3,2}:= \sqrt{-3} (s^2+106
s+2917) (s+35)/144\end{cases}\\
&P_{3,3}:=(x_{3,3} ,y_{3,3} ),~\begin{cases}
&x_{3,3} :=s^2/48 +793/48+35s/24  +  (s+35)\sqrt{-3}/2\\
&y_{3,3} :=
 (-1+\sqrt{-3}) (s+35) (s+6\sqrt{-3}+53)/8\end{cases}\\
&P_{3,4}:=(x_{3,4} ,y_{3,4} ),~\begin{cases} 
&x_{3,4} :=s^2/48  +793/48+35s/24  -   (s+35)\sqrt{-3}/2\\
&y_{3,4} :=
  -  (1+\sqrt{-3}) (s+53-6\sqrt{-3}(s+35)/8\end{cases}
\end{eqnarray}
Thus they generate a subgroup isomorphic to 
$\bfZ/3\bfZ+\bfZ/3\bfZ$. We can take the generators $P_{3,1}, P_{3,2}$
for example.
Thus by 
\cite{Kenku},  $(D_s)_{tor}(\bfQ(\sqrt{-3}))$ is isomorphic to 
  one of the
following.

 (a) $\bfZ/3\bfZ+\bfZ/3\bfZ$, (b) $\bfZ/3\bfZ+\bfZ/6\bfZ$
and (c) $\bfZ/6\bfZ+\bfZ/6\bfZ$. 

\noindent
The case (b) is forgotten in the
list of \cite{Kenku} by an obvious type mistake.
By the same discussion
as in 5.1,  there exists  $P\in D_s$  with order $6$ and 
$2P=P_{3,1}$ if and only if 
\begin{eqnarray*}
&\De(s,m):=s^3+85 s^2-4 m s^2-568 m s+1555 s-16 m^2 s
-1136 m^2\\
&\qquad \qquad -15465-20164 m+64 m^3=0
\end{eqnarray*}
Fortunately the variety $\De=0$ is again rational and 
we can parametrize it as
\begin{eqnarray}
&s=\xi_6(t),\quad\xi_6(t)  := -(27 t^3-1304 t^2+17920
t-71680)/(t-8) (t-16)^2\\ &m=\psi(t),\quad
\psi(t) := -(-128 t^2+3
t^3+1536 t-6144)/(t-8) (t-16)^2
\end{eqnarray}
It turns out that the condition for the existence of $Q\in D_s$
with $2 Q=P_{3,2}$ is the same with the existence of $P,~2P=P_{3,1}$.
Assume that $s=\xi_6(t)$. Then by an easy computation, we get 
$P=(x_{6,1},y_{6,1})$ and $Q=(x_{6,2},y_{6,2})$ where
\begin{eqnarray*} 
&x_{6,1} :=-\frac{1}{3}
\frac{(-3072 t^5+11796480 t^2+86016 t^4-1327104 t^3
-56623104 t+113246208+47 t^6)}{(t-8)^2 (t-16)^4}\\
&y_{6,1} :=
\frac{-4 t^3 (t^2-24 t+192) (7 t^2-144 t+768)}
{(t-16)^5 (t-8)^2}\\
 &x_{6,2}
:=
\frac13\frac{(37 t^6-2016 t^5+40704 t^4-294912 t^3
-1179648 t^2+28311552 t-113246208)}
{(t-8)^2 (t-16)^4}\\
 &y_{6,2 }:=-\frac 87
\frac{\sqrt{-3}(t-12) (t-12-4 \sqrt{-3}) (7 t-72+8 \sqrt{-3})
 (7 t-72-8 \sqrt{-3}) t (t-12+4\sqrt{-3})}
{(t-16)^3 (t-8)^3}
\end{eqnarray*}
It is easy to see by a direct computation that 
$3P=3Q=(\al,0)$ where 
\[\al:=
-\frac 2 3\frac{(t^2-48 t+384) 
(13 t^4-528 t^3+8064 t^2-55296 t+147456)}
{(t-8)^2 (t-16)^4}\]
and $Q-P=P_{3,3}$.
Now we claim that 
\begin{claim}
$(D_s)_{tor}(\bfQ(\sqrt{-3}))=\bfZ/3\bfZ+\bfZ/6\bfZ$ with 
generators $P_{3,3}$ and $P$. 
\end{claim}
In fact, if the torsion is  $\bfZ/6\bfZ+\bfZ/6\bfZ$,
there exist  three elements of order two.
However $f_0(x):=f(x,0)$ factorize as 
$(x-\al) f_{0,0}(x)$ and   their discriminants
are given by 
\begin{eqnarray*}
&\De_x f_0 :=\frac{2048 t^6 (t-12)^3 (t^2-24 t+192)^3
(7 t^2-144 t+768)^6}
{(t-8)^9 (t-16)^{18}}\\
&\De_x f_{0,0} :=165888 (t-12)^3 (t^2-24 t+192)^3 (t-8)^7 (t-16)^8
\end{eqnarray*}
 Consider quartic $Q_4: g(t,v):=165888 (t-12) (t^2-24 t+192) (t-8)-v^2=0$.
Thus $D_s$ has three two torsion elements if and only if 
the quartic $g(t,v)=0$ has  $\bfQ(\sqrt{-3})$-point
$(t_0,v_0)$ with 
$t_0\ne 8,16,12,12\pm 4\sqrt{-3}$. The proof of Claim is reduces to:
\begin{assertion}
There are no such point on $Q_4$.\end{assertion}
{\em Proof.}
By an easy birational change of coordinates, 
$g(t,v)=0$ is equivalent to the elliptic curve
$C:=\{ x^3+1/16777216-y^2=0\}$.
We see that $C$ has two element of order three,
$(0,\pm 1/4096)$  and three two-torsion 
$(-1/256,0), ( 1/512-1/512\sqrt{-3}, 0)$
and $(1/512+1/512\sqrt{-3},0)$. Again by \cite{Kenku},
$C_{tor}(\bfQ(\sqrt{-3}))=\bfZ/2\bfZ+\bfZ/6\bfZ$.
As the rank of $C$ is 0 (\cite{S-Z}), there are exactly 12
points   on $C$.
They correspond to either zeros or poles of 
$\De_x(f_0)$. This implies that the quartic $Q_4$ has no
non-trivial points and thus $C$ does not have three 2-torsion points.
This completes the proof of the  Assertion and thus also proves the
Claim. \qed

Now we formulate our result as follows.
Let $A_6=\{s=\xi_6(t);t\in \bfQ(\sqrt{-3})\}$
and $\Si_6:=\xi_6\inv(\Si)$ is given by
$\Si_6=\{ 8,16,0,12,12\pm 4\sqrt{-3},(72\pm 8\sqrt{-3})/7\}$.
\begin{theorem} \label{quadratic}The Mordell-Weil torsion of $D_s$ is
given by
\[(D_s)_{tor}(\bfQ(\sqrt{-3}))=
\begin{cases}
&\bfZ/3\bfZ+\bfZ/3\bfZ \quad s\in \bfQ(\sqrt{-3})-A_6\cup \Si\\
&\bfZ/6\bfZ+\bfZ/3\bfZ \quad s=\xi_6(t)\in A_6,~t\in \bfQ(\sqrt{-3})-\Si_6
\end{cases}
\]
The j-invariant is given by 
\[j(D_s)=\frac 1{64}
\frac{(s+47)^3 (s+71)^3 (s^2+70 s+1657)^3}
{(s+35)^3 (s^2+106 s+2917)^3}\]
\end{theorem}
\subsection{Examples}

(A) First we consider the case of elliptic curves  $C_s$. In
 the following examples, we give only the values of parameter $s$
as the coefficients are fairly big.
The corresponding 
Weierstrass normal forms are obtained by (\ref{coeff-torus}).

1. $s=54$. The curve $C_{54}$  with torsion group $\bfZ/6\bfZ$  has been studied in 
\S \ref{s=54}.

2. Take $r=3,~s=\vphi_{6,2}(3)=343/9$. Then the torsion group is 
isomorphic to $\bfZ/6\bfZ+\bfZ/2\bfZ$ with generators
$P_2=(-55223/972,-588245/486)$ and $R=(88837/972,0)$.
The j-invariant is given by $7^3\cdot 127^3/2^2\cdot 3^6\cdot 5^2$.

3. Take $t=-3,~s=\vphi_9(-3)=1/216$. 
Then the torsion group is isomorphic  to $\bfZ/9\bfZ$ and 
the generator $P_3=(289/559872, -7/419904)$.
The j-invariant is $71^3\cdot 73^3/2^9\cdot 3^9\cdot 7^3\cdot 17$.

4. Take $\nu=3,~s=\vphi_{12}(3)=-27/80$. Then 
the torsion is isomorphic to $\bfZ/12\bfZ$ with generator
$P_4=(-2997/25600,-6561/102400)$.
The j-invariant is $-11^3\cdot 59^3/2^{12}\cdot 3\cdot 5^3$.

\noindent
(B) We consider elliptic curves $D_s$ defined over $\bfQ(\sqrt{-3})$.
The normal form is given by (\ref{coeff-nontorus}).

5. Take $s=1$. Then  $(D_1)_{tor}(\bfQ(\sqrt{-3}))=\bfZ/3\bfZ+\bfZ/3\bfZ$
and the generators are $(x_{3,1},y_{3,1})=(108,756)$ and
 $(x_{3,2},y_{3,2})=(-144,756\sqrt{-3})$.
The j-invariant is $2^{15}3^3/7^3$.

6. Take $t=4$ and $s=-299/9$. Then the torsion is isomorphic to 
$\bfZ/6\bfZ+\bfZ/3\bfZ$. The generators can be taken as
$(x_{6,1},y_{6,1})=(-2351/243,-532/243)$ and 
$(x_{3,3},y_{3,3})=(8\sqrt{-3}/9-2171/243,-680/81+248\sqrt{-3}/81)$.
The j-invariant is given by \nl
$5^3\cdot 17^3\cdot 31^3\cdot 2203^3/2^6\cdot 3^6\cdot 7^3\cdot 19^6$.

\subsection{Appendix. Parametrization of rational curves}
Parametrizations of a rational curves are always possible
and there exists  even  some programs to find a parametrization
on Maple V. For the detail, see \cite{Ab-Ba} and \cite{vH} for example.
In our case, it is easy to get a parametrization by a direct computation.
For a rational curves with degree less than or equal four is easy.
For other case, we first decrease the degree, using suitable bitational maps.
We give a brief indication. We remark here that the parametrization
is unique up to a linear fractional  change of the parameter.

\noindent
(1) For the parametrization of 
$s^3-32 s^2-2 m s^2-4 m^2 s+8 m^3=0$,
put $m=us$.

\noindent
(2) For the parametrization of 
\begin{eqnarray*}
&R_3(m,s):=512 m^9+768 m^8 s-512 m^6
s^3-1536 m^6 s^2-192 s^4 m^5 \\& -6144 m^5 s^3-6528 m^4 s^4+96 s^5
m^4-12288 m^3 s^4 -2048 m^3 s^5+64 s^6 m^3 
+480 s^6 m^2\\&-15360 s^5 m^2 -6144 s^6 m 
+384 s^7 m-6 s^8 m+56 s^8-512 s^6-768 s^7-s^9 =0
\end{eqnarray*}
put successively $s=s_1/m_1$ and $m=1/m_1$,
then put $n_1=n_2/s_1^2$, then $s_1=s_2-2$ and $n_2=n_4 s_2$. This
changes degree of our curve to be 6. Then $s_2+s_3-4$ and
$n_4=n_5+2$ and $n_5=n_6 s_3$. This changes our curve into 
a quartic. Other computation is easy.

\subsection{Further remark.}
Professor A. Silverberg kindly communicated us 
about the paper \cite{Silverberg}.
He gave a universal family for $\bfZ/3\bfZ+\bfZ/3\bfZ$ over 
$\bfQ(\sqrt{-3})$, which is given by 
$A(u): y^2=x^3+a_0(u) x+b_0(u)$ where
\[
a_0(u)=-27u(8+u^3),\quad b_0(u)=-54(8+20u^3-u^6)\]
and  the  subfamily, given by 
$u=(4+\tau^3)/(3\tau^2)$,
describes elliptic curves with torsion $\bfZ/6\bfZ+\bfZ/3\bfZ$.
 Again by an easy computation, we can show that 
by the change of parameter 
$s=-47+12 u$ we can  identify $D_s$ and $A(u)$. Our subfamily for 
$\bfZ/6\bfZ+\bfZ/3\bfZ$ is also the same with that of \cite{Silverberg}
by the fractional change of parameter:
$t=8(\tau-2)/(\tau-1)$.

\vspace{.2cm}
We would like to thank    H. Tokunaga for 
the valuable discussions and  informations about elliptic
fibrations and also  to  K. Nakamula and  T. Kishi for the information
about elliptic curves over a number field. I am also  gratefull 
to SIMATH for many computations.

\end{document}